\documentclass{amsart}
\usepackage{amsfonts}
\usepackage{amsmath}
\usepackage{amssymb}

\newcommand{\heute}{26 February 2004}

\theoremstyle{plain}
\newtheorem{theorem}{Theorem}[section]
\newtheorem{bigthm}{Theorem}

\newtheorem{lemma}[theorem]{Lemma}

\theoremstyle{remark}

\newtheorem*{rk}{Remark}

\newcommand{\enref}[1]{\textup{(\ref{enum:#1})}}
\newcommand{\dashHsop}{\textup{(\ref{enum:ResHsop}${}'$)}}

\newcommand{\Gro}[1]{Gr\"ob\-ner}

\newcommand{\Ess}{\operatorname{Ess}}
\newcommand{\Ann}{\operatorname{Ann}}

\newcommand{\Id}{\operatorname{Id}}
\newcommand{\eps}{\varepsilon}

\DeclareMathOperator{\Res}{Res}

\newcommand{\bpsi}{\bar{\psi}}

\begin{document}

\title{The essential ideal is a Cohen--Macaulay module}
\author{David J. Green}
\address{Department of Mathematics \\
University of Wuppertal \\ D-42097 Wuppertal \\
Germany}
\email{green@math.uni-wuppertal.de}
\date{\heute}

\begin{abstract}
\noindent
Let $G$ be a finite $p$-group which does not contain a rank two elementary
abelian $p$-group as a direct factor. Then
the ideal of essential classes in the mod-$p$ cohomology ring of~$G$
is a Cohen--Macaulay module whose Krull dimension is the $p$-rank of
the centre of~$G$.
This basically answers in the affirmative a question posed by J.~F.
Carlson (Question~5.4 in~\cite{Carlson:Problems}).
\end{abstract}

\maketitle

\section*{Introduction}
\noindent
Let $G$ be a finite $p$-group and $k$ a field of characteristic~$p$.
By definition, the essential ideal $\Ess(G)$ in the cohomology ring
$H^*(G,k)$ consists of those cohomology classes whose restriction
to every proper subgroup is zero. The essential ideal is involved in the
cohomological characterization of $p$-groups
% with no non-central exponent~$p$ elements~\cite{AdKa:Ess},
whose exponent~$p$ elements are all central~\cite{AdKa:Ess},
and it has featured prominently
in investigations concerning the depth of the graded commutative
ring $H^*(G,k)$ (see~\cite{Carlson:DepthTransfer}).

The essential ideal also plays a key role in Carlson's
method~\cite{Carlson:Tests} for computing $H^*(G,k)$.
In this paper we shall show that the most important case of a question
posed by Carlson (Question~5.4 in \cite{Carlson:Problems})
can be answered using existing techniques. That is, we shall prove
the following theorem:

\begin{bigthm}
\label{ThmA}
Let $k$ be a field of characteristik~$p$,
and let $G$ be a finite $p$-group which does not have the elementary abelian
$p$-group of order~$p^2$ as a direct factor.
If the essential ideal $\Ess(G)$ in $H^*(G,k)$ is nonzero, then
it is a Cohen--Macaulay module with Krull dimension equal to the
$p$-rank of the centre of~$G$.
\end{bigthm}

\noindent
Another formulation of the result is as follows:

\begin{bigthm}
\label{theorem:main}
With the assumptions of Theorem~\ref{ThmA}, denote
by~$C$ the largest central elementary abelian subgroup of~$G$.
Suppose that $\zeta_1,\ldots \zeta_d$ are homogeneous
elements of $H^*(G,k)$ such that the restrictions
$\Res_C(\zeta_1),\ldots,\Res_C(\zeta_d)$ form a homogeneous system
of parameters for $H^*(C,k)$.
Then the essential ideal $\Ess(G)$ in
$H^*(G,k)$ is free and finitely generated as a module over the polynomial
algebra $k[\zeta_1,\ldots,\zeta_d]$.
\end{bigthm}

Carlson's method for computing $H^*(G,k)$ depends upon the essential ideal
being free and finitely generated. Up till now the computer program had
to check this property for each group it was presented with. Now that
the question is a theorem, programs based on Carlson's method
should run faster.
Note that the same question was also mentioned by H.~Mui several years earlier,
but in a much weaker form~\cite{Mui:Essay}.

The equivalence of Theorems \ref{ThmA}~and \ref{theorem:main} is demonstrated
in Lemma~\ref{lemma:tfae}, which also shows that it suffices to
find one sequence $\zeta_1,\ldots,\zeta_d$ satisfying the assumptions
and conclusion of Theorem~\ref{theorem:main}.
Finite generation of $\Ess(G)$ is proved in Lemma~\ref{lemma:fg}
for all finite $p$-groups~$G$; and freeness is proved for the groups
under consideration in~\S\ref{section:free},
using the fact that the cohomology ring of~$G$ is a comodule over the
cohomology ring of any central subgroup.
From now on we shall write $H^*(G)$ for $H^*(G,k)$.

\section{Group cohomology and commutative algebra}
\noindent
We shall demonstrate that Theorems \ref{ThmA}~and \ref{theorem:main}
are equivalent by recalling several well documented facts about group
cohomology and commutative algebra.
The following lemma is certainly not new.

\begin{lemma}
\label{lemma:3ideals}
Let $k$~be a field of characteristic~$p$ and $G$~a finite $p$-group whose
essential ideal $\Ess(G)$ is nonzero.
Denote by $C$ the largest central elementary abelian subgroup of~$G$.
Let $A$, $I$ and $J$ be the following ideals in $H^*(G)$:
\begin{itemize}
\item
$A$~is the annihilator ideal $\Ann_{H^*(G)}(\Ess(G))$;
\item
$I$~is the kernel of restriction $H^*(G) \rightarrow H^*(C)$; and
\item
$J$~is generated by all transfer classes from all proper subgroups.
\end{itemize}
Then $\sqrt{A} = \sqrt{I} = \sqrt{J}$, and this radical ideal is a prime ideal.
\end{lemma}

\begin{proof}
The cohomolgy ring of an elementary abelian $p$-group is a polynomial
algebra ($p=2$) or the tensor product over~$k$ of a polynomial algebra with
an exterior algebra ($p$~odd). So $H^*(C)/\sqrt{0}$ is a polynomial algebra,
and $\sqrt{I}$ a prime ideal. Quillen proved that $x \in H^*(G)$ is nilpotent
if its restriction to every elementary abelian subgroup is nilpotent
(Coroll.~8.3.4 in~\cite{Evens:book})\@. So $\sqrt{I}$ consists of those
classes with nilpotent restriction to the centre~$Z(G)$. It follows from
Carlson's theorem on varieties and transfer (\S10.2 in~\cite{Evens:book})
that $\sqrt{I}=\sqrt{J}$.

Let $H \leq G$ be a subgroup. Then restriction makes $H^*(H)$ a module
over $H^*(G)$, and transfer from $H^*(H)$ to $H^*(G)$ is a homomorphism
of $H^*(G)$-modules. Hence transfer classes and essential classes annihilate
each other, which means that $J \subseteq A$.

To finish the proof we will use an argument of Broto and
Henn~\cite{BrHe:Central} to show that every element
of $H^*(G) \setminus \sqrt{I}$ is regular, from which it follows
that $A \subseteq \sqrt{I}$.
Group multiplication $\mu \colon G \times C \rightarrow G$
is a group homomorphism, and so induces a map of $k$-algebras
$\mu^* \colon H^*(G) \rightarrow H^*(G) \otimes_k H^*(C)$.
Moreover, $(\Id_G \times \eps_C) \circ \mu^* = \Id_G$ and
$(\eps_G \otimes \Id_C) \circ \mu^* = \Res^G_C$,
where $\Id_G \colon H^*(G) \rightarrow H^*(G)$ is the identity
map and $\eps_G \colon H^*(G) \rightarrow k$ the augmentation.
So $\mu^*$ is a split monomorphism; and $\mu^*(x)
\in 1 \otimes \Res_C(x) + H^{>0}(G) \otimes H^*(C)$ for all $x \in H^*(G)$.
So if $\Res_C(x)$ is non-nilpotent, then it is a regular element of
$H^*(C)$, and it follows that $\mu^*(x)$ is a regular element of
$H^*(G) \otimes_k H^*(C)$. As $\mu^*$ is injective, it follows that
$x \in H^*(G)$ is regular.
\end{proof}

\begin{lemma}
\label{lemma:fg}
Let $G$ be a finite $p$-group, and $k$ a field of characteristic~$p$.
Denote by $C$ the largest central elementary abelian subgroup of~$G$.
Suppose that $\zeta_1,\ldots \zeta_d$ are homogeneous
elements of $H^*(G)$ whose restrictions to~$C$
form a homogeneous system of parameters for $H^*(C)$.
Then $\Ess(G)$ is finitely generated as a module over the polynomial
algebra $R = k[\zeta_1,\ldots,\zeta_d]$.
\end{lemma}

\begin{proof}
The Evens--Venkov theorem (Coroll.\@ 7.4.6 in~\cite{Evens:book}) states that
$H^*(G)$ is a finitely generated $k$-algebra.
Let $A,I,J$ be as in Lemma~\ref{lemma:3ideals}, and recall from the proof
of that lemma that $J \subseteq A$.
So $\Ess(G)$ is a module over $H^*(G)/J$; and it is finitely generated since
$H^*(G)$ is noetherian.
So it would suffice to show that $H^*(G)/J$ is a finitely generated $R$-module.

If $S \subseteq T$ are finitely generated (graded-)commutative $k$-algebras
and $L \subseteq T$ an ideal such that $T/\sqrt{L}$ is a finitely generated
$S$-module, then $T/L$ is a finitely generated $S$-module too. So it suffices
to show that $H^*(G)/\sqrt{J}$ is a finitely generated $R$-module.
But we have seen that $\sqrt{J}=\sqrt{I}$;
and by assumption $H^*(G)/I$ is a submodule
of the noetherian $R$-module $H^*(C)$.
\end{proof}

\begin{lemma}
\label{lemma:tfae}
Let $G$ be a finite $p$-group and $k$ a field of characteristic~$p$.
Denote by $C$ the largest central elementary abelian subgroup of~$G$.
Assuming that $\Ess(G)$ is nonzero, the following three statements
are equivalent:
\begin{enumerate}
\item
\label{enum:CM}
$\Ess(G)$ is a Cohen--Macaulay $H^*(G)$-module.
\item
\label{enum:exists}
There is a sequence $\zeta_1,\ldots \zeta_d$ of homogeneous
elements of $H^*(G)$ such that
\begin{enumerate}
\item
\label{enum:ResHsop}
the restrictions $\Res_C(\zeta_1),\ldots,\Res_C(\zeta_d)$ form a homogeneous
system of parameters for $H^*(C)$; and
\item
\label{enum:FgFree}
the essential ideal $\Ess(G)$ is free and finitely generated as a module
over the polynomial algebra $k[\zeta_1,\ldots,\zeta_d]$.
\end{enumerate}
\item
\label{enum:forall}
Every sequence $\zeta_1,\ldots \zeta_d$ which satisfies
condition~\enref{ResHsop}
also satisfies condition~\enref{FgFree}.
\end{enumerate}
If these equivalent conditions hold, then $d$ and the Krull dimension
of the module $\Ess(G)$ both coincide with the $p$-rank of the centre of~$G$.
\end{lemma}

\begin{proof}
Let $A,I,J$ be as in Lemma~\ref{lemma:3ideals}\@.
First we shall see that condition~\enref{ResHsop} is equivalent to the
following condition:
\begin{enumerate}
\item[\dashHsop]
The images of $\zeta_1,\ldots,\zeta_d$ in $H^*(G)/A$
are algebraically independent over~$k$, and $\Ess(G)$ is a finitely
generated module over $R=k[\zeta_1,\ldots,\zeta_r]$.
\end{enumerate}
Since $\sqrt{I}=\sqrt{A}$, the images of $\zeta_1,\ldots,\zeta_d$ are
algebraically independent in
$H^*(G)/A$
if and only if they are algebraically independent in
$H^*(G)/I \subseteq H^*(C)$.
Lemma~\ref{lemma:fg}
shows that \enref{ResHsop}~implies that $\Ess(G)$ is a finitely generated
$R$-module. Conversely, \dashHsop~implies that $\Ess(G)$ and therefore
$H^*(G)/A$ are finitely generated $R$-modules. Since $\sqrt{I} = \sqrt{A}$,
it follows that $H^*(G)/\sqrt{I}$ is a finitely generated $R$-module.
So just as in the proof of Lemma~\ref{lemma:fg}, $H^*(G)/I$ is a
finitely generated $R$-module.
But recall from Corollary 7.4.7 of~\cite{Evens:book} that
$H^*(C)$ is finite over the image $H^*(G)/I$ of restriction from~$H^*(G)$.
So $H^*(C)$ is a finitely generated $R$-module.
Hence \enref{ResHsop}~and \dashHsop{} are indeed equivalent.

But the characterisation of Cohen--Macaulay modules
in Theorem~4.3.5 of~\cite{Benson:PolyInvts} states precisely that~\enref{CM},
\enref{exists}~and \enref{forall} are equivalent if one replaces
\enref{ResHsop}~by \dashHsop{} in \enref{exists}~and \enref{forall}\@.
Moreover, $d=\dim(\Ess(G))$ would
follow too. As $H^*(C)$ is polynomial for $p=2$ and polynomial tensor
exterior for $p$~odd, and the dimension is equal to the rank of~$C$ in
both cases, $d$ is the $p$-rank of the centre of~$G$.
\end{proof}

\begin{rk}
Assuming that the equivalent conditions \enref{ResHsop},~\dashHsop{}
hold, apply the same characterisation of Cohen--Macaulay modules to $H^*(C)$.
Recalling that $H^*(C)$ is a polynomial algebra for $p=2$, and polynomial
tensor exterior for $p$~odd,
one deduces that $H^*(C)$ is a free and finitely generated $R$-module.
\end{rk}

\section{Freeness}
\label{section:free}
\noindent
As in~\cite{BrHe:Central},
observe that group multiplication $\mu \colon G \times C \rightarrow G$
turns $H^*(G)$ into a comodule over the coalgebra $H^*(C)$.
As was noted above, the comodule structure map $\mu^* \colon H^*(G) \rightarrow
H^*(G) \otimes_k H^*(C)$ is simultaneously a map of $k$-algebras.

The following lemma is the key to this paper.
It is perhaps the natural level of generality
for a result that has been known in various special cases for
some time (Lemma~3.1 in~\cite{AdMi:Central}, Proposition~5.2
in~\cite{Carlson:Problems} and Lemma~1.2 in~\cite{depth})\@.
The lemma bears a striking resemblance to a classical fact about Hopf
algebras (Theorem~4.1.1 in~\cite{Sweedler}),
but I do not know see to derive the the former as a corollary of
the latter.

\begin{lemma}
\label{lemma:subcomodFree}
Let $G,C,R$ be as in Lemma~\ref{lemma:fg}, and
let $F = \bigoplus_{i \geq i_0} F_i$ be the free graded $H^*(G)$-module 
on generators $e_1,\ldots,e_s$
(not necessarily in degree zero). Hence $F$~is a comodule over~$H^*(C)$,
with structure map~$\psi$ given by
$\psi (ae_j) = \mu^*(a) \cdot (e_j \otimes 1) \in F \otimes_k H^*(C)$.
Let $N \subseteq M$ be graded submodules of~$F$ which are simultaneously
subcomodules. Then $M/N$ is a free $R$-module.
\end{lemma}

\begin{proof}
Define $\bpsi \colon M/N \rightarrow M/N \otimes_k H^*(C)$
by $\bpsi(a + N) = \psi(a) + N \otimes_k H^*(C)$. By the assumptions,
$(M/N, \bpsi)$ is an $H^*(C)$-comodule. Moreover $\bpsi$ is a map
of $H^*(G)$-modules if one gives $M/N \otimes_k H^*(C)$ the $H^*(G)$-module
structure induced from the obvious $H^*(G) \otimes_k H^*(C)$-module
structure by~$\mu^*$. In fact, $\bpsi$ is a split monomorphism
of $H^*(G)$-modules, the splitting map being $\Id_{M/N} \otimes \varepsilon$,
where $\varepsilon \colon H^*(C) \rightarrow k$ is the augmentation map.

Inclusion of $R$~in $H^*(G)$ induces an $R$-module structure on
$M/N \otimes_k H^*(C)$. We shall show that this $R$-module is free.
It then follows that $\zeta_1,\ldots,\zeta_z$ is a regular sequence for
this $R$-module, and therefore for the $R$-module $M/N$~too.
Hence $M/N$~is a free $R$-module.

For $i \geq i_0$ set $S_i := \bigoplus_{j \geq i} (M/N)_i$,
and let $T_i$ be the $R$-submodule $S_i \otimes_k H^*(C)$
of $M/N \otimes H^*(C)$. Then
$M/N \otimes_k H^*(C)
= T_{i_0} \supseteq T_{i_0+1} \supseteq T_{i_0+2} \supseteq \cdots$
and $\bigcap_{i \geq i_0} T_i = \{0\}$.
Now, $T_i/T_{i+1}$ is $H^*(C)$-module, and as such it is free of rank
$\dim_k (M/N)_i$. Moreover the $R$-module structure is induced from
this $H^*(C)$-module structure by inclusion $R \hookrightarrow H^*(G)$
followed by restriction to $H^*(C)$, since for all
$\rho \in R$, $x \in S_i$ and $y \in H^*(C)$ one has
$\rho \cdot (x \otimes y + T_{i+1})
= x \otimes \Res_C(\rho) \cdot y + T_{i+1}$.
We observed above that $H^*(C)$~is a free $R$-module, and so
$T_i/T_{i+1}$ is a free $R$-module. By induction,
$T_{i_0} / T_i$ is a free $R$-module for all~$i$. Since $T_i$ is confined to
degree${} \geq i$, this means that $T_{i_0} = M/N \otimes_k H^*(C)$ is
itself a free $R$-module.
\end{proof}

\begin{lemma}
\label{lemma:seqsExist}
In Lemma~\ref{lemma:tfae} there is always a sequence
$\zeta_1,\ldots,\zeta_d$ satisfying condition~\enref{ResHsop}.
\end{lemma}

\begin{proof}
Restriction turns $H^*(C)$ into an $H^*(G)$-module;
recall from Corollary 7.4.7 of~\cite{Evens:book} that
this module is finitely generated.
The result follows by the graded case of Noether Normalization
(Theorem~2.2.7 in~\cite{Benson:PolyInvts}) applied to the $H^*(G)$-module
$H^*(C)$: recall that the annihilator ideal of this $H^*(G)$-module
is the kernel of the restriction map.

One explicit sequence is $\zeta_i = m$th Chern class of the regular
representation of~$G$, for $m=2(p^n-p^{n-i})$, where $|G|=p^n$.
This construction features in Venkov's topological proof of the Evens--Venkov
theorem (\S3.10 in~\cite{Benson:II}).
\end{proof}

\begin{proof}[Proof of Theorems \ref{ThmA}~and \ref{theorem:main}]
Lemma~\ref{lemma:tfae} shows that the two theorems are equivalent,
and that to prove them it suffices to prove the existence of a sequence
$\zeta_1,\ldots,\zeta_d$ satisfying conditions \enref{ResHsop}~and
\enref{FgFree}\@. Lemma~\ref{lemma:seqsExist} demonstrates that
there is always a sequence satisfying~\enref{ResHsop},
and by Lemma~\ref{lemma:fg} the finite generation part of~\enref{FgFree}
is an automatic consequence of~\enref{ResHsop}.
We shall adopt the notation of Lemma~\ref{lemma:fg}\@.

Assume first that the cyclic group of order~$p$ is not a direct factor of~$G$.
Then every maximal subgroup of~$G$ contains~$C$, and so
$\Ess(G)$ is a subcomodule of the $H^*(C)$-comodule $H^*(G)$.
Hence $\Ess(G)$ is a free $R$-module by Lemma~\ref{lemma:subcomodFree},
so the theorems are proved for this~$G$.

Now suppose that $G = H \times C_p$, and that $H$ has no direct factor which
is cyclic of order~$p$.
Denote by $C'$ the largest central elementary abelian subgroup of~$H$.
Then $C = C' \times C_p$, and $\Ess(H)$ is a $H^*(C')$-comodule by
the above argument. Applying Lemma~\ref{lemma:seqsExist} to $H,C'$
one obtains a sequence $\zeta_1,\ldots,\zeta_{d-1} \in H^*(H)$ which satisfies
condition~\enref{ResHsop} of Lemma~\ref{lemma:tfae} for $H,C'$.
Let $t$ be a homogeneous system of parameters for the one-dimensional
$k$-algebra $H^*(C_p)$.
Then $\zeta_1,\ldots,\zeta_{d-1},t$ satisfies condition~\enref{ResHsop}
for~$G,C$.
Set $S = k[\zeta_1,\ldots,\zeta_{d-1}] \subseteq H^*(H)$
and $R = S[t] \subseteq H^*(G)$.
Observe that $\Ess(H)$ is a free $S$-module by the case we have already
proved.
Let $b_1,\ldots,b_n$ be a finite free generating set of the
$k[t]$-module $H^*(C_p)$. We take each~$b_{\ell}$ to be homogeneous.
Let $B \subseteq H^*(C_p)$ be the $k$-vector space with basis
$b_1,\ldots,b_n$. Then $H^*(C_p) = \bigoplus_{j \geq 0} t^j B$.

Set $F_i := H^*(H) \otimes_k \bigoplus_{j=0}^i t^j B$.
Define $M_i$~and $N_i$ by $M_i := \Ess(G) \cap F_i$ and
$N_i := M_{i-1} + M_{i-1} t$.
Both $M_i$~and $N_i$ are submodules
of the free and finitely generated $H^*(H)$-module $F_i$. Moreover they
are subcomodules of the $H^*(C')$-comodule $F_i$,
as $C'$~lies in every maximal subgroup of~$G$.
So $M_i/N_i$ is a free $S$-module by Lemma~\ref{lemma:subcomodFree}.
Observe that the $R$-modules $k[t]N_i$ and $k[t]M_{i-1}$ coincide.
Moreover, the noetherian $R$-module $\Ess(G)$ is the union
of the $R$-modules $k[t]M_i$, and so $k[t]M_i=\Ess(G)$
for large enough~$i$.

It would suffice to show that each $k[t]M_i$ is a free $R$-module,
and this follows by induction if we can show that each
$(k[t] M_i)/(k[t] N_i)$ is a free $R$-module.
(Note that $N_0 = \{0\}$.)
We shall show this by proving that $(k[t] M_i)/(k[t] N_i)$ is the
free $S[t]$-module $k[t] \otimes_k (M_i/N_i)$.

Note that $\Ess(G) \subseteq \Ess(H) \otimes_k H^*(C_p)$,
since $K \times C_p$ is a maximal subgroup of $G$ for
each maximal subgroup $K$~of $H$. (For $H=1$ one has
$\Ess(1)=k=H^*(1)$.)

Let $y$~be an element of $M_i$ in degree~$m$. Then
there are homogeneous classes $a_{j\ell} \in \Ess(H)$
for $0 \leq j \leq i$ and $1 \leq \ell \leq n$ with
$y = \sum_{j=0}^i \sum_{\ell=1}^n a_{j\ell} \otimes b_{\ell} t^j$.
Call $\sum_{\ell} a_{i\ell} \otimes b_{\ell}$ the leading
coefficient of~$y$, and let $I_i \subseteq \Ess(H) \otimes_k B$
denote the set of all such leading
coefficients.
Then $I_i$ is an $H^*(H)$-module, and $y$~lies
in $N_i$ if and only if its leading coefficient lies in $I_{i-1}$.
So we may pick a basis for a complement of the subspace $N_i$~of $M_i$
consisting of classes whose leading terms constitute a basis for
a complement of the subspace $I_{i-1}$~of $I_i$. So the natural surjective
map of $S[t]$-modules $k[t] \otimes_k (M_i/N_i) \rightarrow
(k[t] M_i)/(k[t] N_i)$ is injective too. But the domain is a free $S[t]$-module.
\end{proof}

\section{Concluding remarks}

\subsection*{Application to computer calculations}
In~\cite{Carlson:Tests}, Carlson describes a series of tests on a
partial presentation for $H^*(G)$ and proves that the presentation is
complete if it passes the tests. This is of crucial importance for the
computer calculation of group cohomology via minimal resolutions, for
otherwise one would never know when to stop.

However Carlson's tests depend on two conjectures about the structure
of the cohomology ring, which have to be checked for the group in
question as part of the calculation.
This means that there could conceivably be some groups where
there is no complete presentation that the tests can detect as being complete.
One of these conjectures concerns the Koszul complex associated
to a homogeneous system of parameters.
The other is now proven, as our Theorem~\ref{theorem:main}: note that Carlson's
tests assume that the cohomology rings of all subgroups are known, and
so the cohomology ring of a product group is known by the K\"unneth theorem.

\providecommand{\bysame}{\leavevmode\hbox to3em{\hrulefill}\thinspace}

\end{document}